\documentclass[11 pt]{amsart}
\usepackage{amsmath,amssymb,amsthm}%\textheight20cm
\newtheorem{Theorem}{Theorem}[section]

\newtheorem{Corollary}[Theorem]{Corollary}

\newtheorem{Lemma}[Theorem]{Lemma}

\theoremstyle{definition}
\newtheorem{Example}[Theorem]{Example}
\newtheorem{Remark}[Theorem]{Remark}

\begin{document}   % End of preamble and beginning of the text

%\doublespacing

\title[Relative Flux Homomorphism]
{Relative Flux Homomorphism in Symplectic Geometry}

\author{Y{\i}ld{\i}ray Ozan}
\address{Department of Mathematics, Middle East Technical University,
06531 \newline Ankara, TURKEY} \email{ozan@metu.edu.tr}
\date{\today}
%\thanks{The author is partially supported by.}

\pagenumbering{arabic}

\begin{abstract}
In this work we define a relative version of the flux
homomorphism, introduced by Calabi in 1969, for a symplectic
manifold.  We use it to study (the universal cover of) the group
of symplectomorphisms of a symplectic manifold leaving a
Lagrangian submanifold invariant. We also show that some quotients
of the universal covering of the group of symplectomorphisms are
stable under symplectic reduction.
\end{abstract}

\maketitle

\section{Introduction}
In 1969 Calabi introduced the flux homomorphism for a symplectic
manifold $(M,\omega)$ (\cite{Cal}). Let
$\textmd{Symp}_0(M,\omega)$ denote the path component of the Lie
group of symplectomorphisms of $(M,\omega)$ containing the
identity. Also, let $\widetilde{\textmd{Symp}}_0(M,\omega)$ denote
the universal covering space of $\textmd{Symp}_0(M,\omega)$.  Note
that an element of $\widetilde{\textmd{Symp}}_0(M,\omega)$ is
represented by a homotopy class of a smooth path $\psi_t$ in
$\textmd{Symp}_0(M,\omega)$ connecting the identity $\psi_0=id$ to
$\psi=\psi_1$, where the homotopies fix the end points of the
paths.  We will denote the homotopy class of $\psi_t$ by $\{
\psi_t \} \in \widetilde{\textmd{Symp}}_0(M,\omega)$.  Any such
$\psi_t$ determines a smooth family of vector fields
$X_t:M\rightarrow T_*M$ defined by the formula
$$\frac{d}{dt}\psi_t=X_t\circ \psi_t.$$  Now the flux homomorphism
$\textmd{Flux}:\widetilde{\textmd{Symp}}_0(M,\omega)\rightarrow
H^1(M,{\mathbb R})$ is defined, for any $\{\psi_t\} \in
\widetilde{\textmd{Symp}}_0(M,\omega)$, as
$$\textmd{Flux}(\{ \psi_t\})=\int_0^1[\imath(X_t)\omega] \ dt.$$
For connected $M$ identifying $H^1(M,{\mathbb R})$ with
$\textmd{Hom}(\pi_1(M),{\mathbb R})$ the cohomology class
$\textmd{Flux}(\{ \psi_t\})$ corresponds to the homomorphism
$\pi_1(M)\rightarrow {\mathbb R}$ defined by
$$\gamma \mapsto \int_0^1 \int_0^1 \omega(X_t(\gamma(s),
\dot{\gamma}(s))) \ ds \ dt$$ for any smooth loop $\gamma
:S^1\rightarrow M$.

In the next section, we will define a relative flux homomorphism
in the presence of a Lagrangian submanifold $L\subseteq M$, mainly
following both the treatment and the notation of Section $10$ of
\cite{McS}.  Most of the results and proofs of this section are
analogous to those of \cite{McS} and we will omit the proofs of
some of them unless the proof includes some new ingredients or
makes the text more comprehensible. In the third section, we will
show that some quotients of
$\widetilde{\textmd{Symp}}_0(M,\omega)$ are stable under
symplectic reduction.

\section{Relative Flux Homomorphism} Throughout this section,
unless stated otherwise, $(M,\omega)$ is a connected closed
symplectic manifold and $L$ a closed Lagrangian submanifold.

Let $\textmd{Symp}(M,L,\omega)$ denote the subgroup of
$\textmd{Symp}(M,\omega)$ consisting of symplectomorphisms leaving
the Lagrangian submanifold $L$ invariant and
$\textmd{Symp}_0(M,L,\omega)$ the path component of
$\textmd{Symp}(M,L,\omega)$ containing the identity.

\begin{Remark} \label{Rem-rel}
In general,
$\textmd{Symp}(M,L,\omega)\cap\textmd{Symp}_0(M,\omega)$ is not
path connected.  Let $S^2$ be the unit sphere in ${\mathbb R}^3$
equipped with the standard symplectic form and $S^1$ be the
intersection of $S^2$ with the $xy-$plane.  Then $S^1$ is a
lagrangian submanifold $S^2$.  Let $\psi_t$, $t\in [0,1]$, denote
the rotation of ${\mathbb R}^3$ about the $x-$axis $t\pi$ radians
and $\psi=\psi_1$.  Since $\psi_{|S^1}$ is orientation reversing
$\psi$ is not in $\textmd{Symp}_0(M,L,\omega)$, even though it
lies in $\textmd{Symp}(M,L,\omega)\cap\textmd{Symp}_0(M,\omega)$
trivially. Hence,
$\textmd{Symp}(M,L,\omega)\cap\textmd{Symp}_0(M,\omega)$ is not
path connected.
\end{Remark}

Identifying the homotopy class of the constant loop at the
identity element of the group ${\textmd{Symp}}_0(M,L,\omega)$ with
that of ${\textmd{Symp}}_0(M,\omega)$, we obtain a canonical map
from $\widetilde{\textmd{Symp}}_0(M,L,\omega)$ onto a subgroup of
$\widetilde{\textmd{Symp}}_0(M,\omega)$, the subgroup of homotopy
classes $\{ \psi_t \}$ in $\widetilde{\textmd{Symp}}_0(M,\omega)$,
where each $\psi_t$ leaves $L$ invariant. Then we have the
following result.

\begin{Lemma}\label{lem-defn}
The flux homomorphism is well defined on $\widetilde{\textmd{\em
Symp}}_0(M,L,\omega)$ and takes values in the kernel of the
restriction map $H^1(M,{\mathbb R})\rightarrow H^1(L,{\mathbb
R})$. If $L$ is connected the flux homomorphism takes values in
$H^1(M,L,{\mathbb R})$.  Moreover, in this case the flux
homomorphism
$$\textmd{\em Flux}:\widetilde{\textmd{\em Symp}}_0(M,L,\omega)\rightarrow
H^1(M,L,{\mathbb R})$$ is onto.
\end{Lemma}

\begin{proof}
The Flux homomorphism is well defined on
$\widetilde{\textmd{Symp}}_0(M,\omega)$ since homotopic loops in
$\textmd{Symp}_0(M,\omega)$ have the same value under the flux
homomorphism (cf. see Lemma $10.7$ of \cite{McS}).  Therefore,
homotopic loops in $\textmd{Symp}_0(M,L,\omega)$ have the same
value and hence the flux homomorphism is well defined on
$\widetilde{\textmd{Symp}}_0(M,L,\omega)$.

Let $\{\psi_t\} \in \widetilde{\textmd{Symp}}_0(M,L,\omega)$ and
$X_t$ the vector field defined by $$\displaystyle
\frac{d}{dt}\psi_t=X_t\circ \psi_t.$$ Since $\psi_t$ leaves $L$
invariant, for any $p\in L$ we have $X_t(p)\in T_pL$. Note that to
prove the first assertion it suffices to show that for any smooth
loop $\gamma :S^1\rightarrow L$ the integral  $$\displaystyle
\int_0^1 \int_0^1 \omega(X_t(\gamma(s), \dot{\gamma}(s))) \ ds \
dt=0.$$  However, this trivially holds since both \
$X_t(\gamma(s))$ and $\dot{\gamma}(s)$ lie in $T_{\gamma(s)}L$,
and $L$ is a Lagrangian submanifold.

Now suppose that $L$ is connected.  Then the exact sequence
$$\rightarrow H^0(M,{\mathbb R})\stackrel{\simeq}{\rightarrow}
H^0(L,{\mathbb R})\stackrel{0}{\rightarrow} H^1(M,L,{\mathbb
R})\rightarrow H^1(M,{\mathbb R})\rightarrow H^1(L,{\mathbb
R})\rightarrow$$ implies that the flux homomorphism takes values
in the isomorphic image of the relative cohomology group
$H^1(M,L,{\mathbb R})$ in $H^1(M,{\mathbb R})$.

To show surjectivity let $u\in \Omega^1(M)$ be a closed $1$-form
which is exact on $L$.  Then $u_{|L}=dh$ for some smooth function
$h:L\rightarrow {\mathbb R}$.  We can extend $h$ first to a
tubular neighborhood and then to whole $M$.  Replacing $u$ by
$u-dh$ we can assume that $u_{|L}=0$.  The closed $1$-form $u$
gives a symplectic vector field $X_t$ on $M$, which is indeed
constant in time. If $\psi_t$ denotes the smooth $1$-parameter
family of symplectomorphisms generated by $X_t$, then clearly
$$\textmd{Flux}(\{ \psi_t \})=\int_0^1[\imath(X_t)\omega] \
dt=\int_0^1 [u] \ dt =[u].$$ We need to show that $\psi_t$ leaves
$L$ invariant. Let $p \in L$ and $Y\in T_pL$ be any vector.  Then
by the construction of $u$ we have $0=u(p)(Y)=\imath(X_t)\omega(p)
(Y)=\omega (X_t(p),Y)$ for any $Y\in T_pL$, which implies that
$X_t(p)\in T_pL$, because $L$ is a Lagrangian submanifold.  This
finishes the proof.
\end{proof}

The following theorem is the relative version of Theorem $10.12$
of \cite{McS}, which describes the kernel of the flux
homomorphism.  First we define the Hamiltonian symplectomorphisms
of $M$ leaving $L$ invariant, denoted $\textmd{Ham}(M,L)$, as the
subgroup of $\textmd{Symp}(M,L)$ consisting of symplectomorphisms
$\psi$ such that there is a Hamiltonian isotopy
$\psi_t:(M,L)\rightarrow (M,L)$, $t\in [0,1]$, such that
$\psi_0=id$ and $\psi_1=\psi$; i.e., $\psi_t$ is a Hamiltonian
isotopy of $M$ such that $\psi_t(L)=L$, for any $t\in [0,1]$.

Since $L$ is Lagrangian a Hamiltonian $H_t:M\rightarrow {\mathbb
R}$ generating the above isotopy is locally constant on $L$.

\begin{Theorem}\label{thm-kernel}
Let $\psi \in \textmd{\em Symp}_0(M,L,\omega)$. Then $\psi$ is a
Hamiltonian symplectomorphism if and only if there exists a
symplectic isotopy, $\psi_t$, $[0,1]\rightarrow \textmd{\em
Symp}_0(M,L,\omega)$ such that $\psi_0=id$, $\psi_1=\psi$ and
$\textmd{\em Flux}(\{ \psi_t\})=0$.

Moreover, if $\textmd{\em Flux}(\{ \psi_t\})=0$ then $\{ \psi_t\}$
is isotopic with fixed end points to a Hamiltonian isotopy through
paths in ${\textmd{\em Symp}}_0(M,L,\omega)$.
\end{Theorem}

The proof of the above theorem is almost the same as that of
Theorem $10.12$ of \cite{McS}, where one has to observe that the
symplectomorphisms and isotopies of $M$ involved in the proof
leave $L$ invariant.

In \cite{McS} the authors define the subgroup $$\Gamma(M)
=\textmd{Flux}(\pi_1(\textmd{Symp}_0(M)))\subseteq H^1(M,{\mathbb
R}).$$  Similarly, we define its relative version $$\Gamma(M,L)
=\textmd{Flux}(\pi_1(\textmd{Symp}_0(M,L,\omega)))\subseteq
H^1(M,L,{\mathbb R})$$ as the image of the fundamental group of
the identity component of the relative symplectic group
$\textmd{Symp}_0(M,L,\omega)$ under the flux homomorphism.  As
$\Gamma(M)$ the group $\Gamma(M,L)$ is a subgroup of
$H^1(M,P_{\omega})$, where $P_{\omega}$ is the  additive subgroup
$[ \omega ] \cdot H_2(M,{\mathbb Z})\subseteq {\mathbb R}$, and
hence both $\Gamma(M)$ and $\Gamma(M,L)$ are countable.

To state the next lemma we need to recall the correspondence
between symplectomorphisms of $(M,\omega)$, which are $C^1-$close
to the identity and closed $1$-forms on $M$, which are close to
the zero form. Consider the symplectic manifold $(M\times M, (-
\omega)\oplus \omega)$ whose diagonal $\Delta \subset M\times M$
is a Lagrangian submanifold, diffeomorphic to $M$.  Then the
Lagrangian neighborhood theorem implies that there is a
symplectomorphism $\Psi :\mathcal{N}(\Delta)\rightarrow
\mathcal{N}(M_0)$ between the neighborhoods of the Lagrangian
submanifolds $\Delta \subset M\times M$ and $M_0\subseteq
(T^*M,\omega_{can})$, where $M_0$ is the zero section of the
$T^*M$ equipped with the canonical symplectic form $\omega_{can}$,
satisfying $\Psi^*(\omega_{can})=(- \omega)\oplus \omega$ and
$\Psi(q,q)=q$, for any $q\in M$.  Then the correspondence
$$(\psi:M\rightarrow M) \mapsto\sigma = \mathcal{C}(\psi)\in
\Omega^1(M)$$ is defined by
$$\Psi(\textmd{graph}(\psi))=\textmd{graph}(\sigma)$$ provided
that $\psi \in \textmd{Symp}(M)$ is sufficiently $C^1-$close to
the identity.

Studying the proof of the existence of a such symplectomorphism
$\Psi$, one realizes that $\Psi$ can be chosen so that
symplectomorphisms of $M$, $C^1-$close to the identity, leaving a
Lagrangian submanifold $L$ invariant correspond to closed
$1$-forms which evaluates zero on $TL$ (cf. see Theorem 3.32 and
Theorem 3.14 of \cite{McS}). Indeed, for the proof of the below
relative version one has to note that a symplectomorphism,
$C^1-$close to the identity, leaves $L$ invariant if and only if
the various $1$-forms involved in the proofs of Theorem 3.32 and
Theorem 3.14 of \cite{McS} vanish on $(L\times L)\cap\Delta$ and
the vector fields are parallel to $L\times L$.

\begin{Lemma}\label{lem-diagonal}
Assume the above notation.  Then there is a symplectomorphism
$\Psi :\mathcal{N}(\Delta)\rightarrow \mathcal{N}(M_0)$ satisfying
$\Psi(q,q)=q$ for any $q\in M$, such that if
$\psi:(M,\omega)\rightarrow (M,\omega)$ is a symplectomorphism,
which is sufficiently $C^1-$close to the identity map, then $\psi$
leaves a Lagrangian submanifold $L$ invariant if and only if
$\sigma_{|T_qL}=0$, for any $q \in L$.
\end{Lemma}

\begin{Remark}\label{rem-diagonal}
The above lemma can be restated as follows: Let $L$ be a
Lagrangian submanifold of a symplectic manifold $(M,\omega)$. Then
there is a symplectomorphism  $\Psi
:\mathcal{N}(\Delta)\rightarrow \mathcal{N}(M_0)$ satisfying
$\Psi(q,q)=q$ for any $q\in M$, which maps a neighborhood of
$(L\times L)\cap \Delta$ in $L\times L$ onto a neighborhood of
$L_0=L\cap M_0$ in the conormal bundle to $L_0\subseteq T^*M$,
which is $\{v\in T^*_qM \ | \ q \in L, \ ~ v_{|T_qL}=0 \}$.
\end{Remark}

An immediate consequence of the above lemma is the following
corollary.

\begin{Corollary}\label{cor-diagonal}
Let $\psi:(M,\omega)\rightarrow (M,\omega)$ and $\sigma =
\mathcal{C}(\psi)\in \Omega^1(M)$ be as in the above lemma. Assume
that $\psi$ leaves a Lagrangian submanifold $L$ invariant. Then
for any $t\in [0,1]$, the symplectomorphism $\psi_t$,
corresponding to the closed $1$-form $t\sigma$, leaves $L$
invariant.
\end{Corollary}

Now we can prove the following relative version of Lemma $10.16$
of \cite{McS}.

\begin{Lemma}\label{lem-gamma}
If $\psi \in \textmd{\em Symp}_0(M,L,\omega)$ is sufficiently
$C^1$-close to the identity and $\sigma= \mathcal{C}(\psi_t)\in
\Omega^1(M)$, then $\psi \in \textmd{\em Ham}(M,L)$ if and only if
$[\sigma]\in \Gamma (M,L)$.
\end{Lemma}

\begin{proof}
Let $\psi_t$, $t \in [0,1]$, be the symplectic isotopy with
$\mathcal{C}(\psi_t)=t\sigma$.  Note that by the above corollary
$\psi_t \in \textmd{Symp}_0(M,L,\omega)$. Lemma $10.15$ of
\cite{McS} implies that $\textmd{Flux}(\{\psi_t\})=-[\sigma]$.
Since $\psi$ is a Hamiltonian symplectomorphism, the path $\psi_t$
extends to a loop $[0,2]\mapsto \textmd{Symp}_0(M,L,\omega)$,
which is generated by Hamiltonian vector fields with end points
$\psi_0=id=\psi_2$.  Moreover, $$\textmd{Flux}(\{\psi_t\}_{0\leq t
\leq 2})=\textmd{Flux}(\{\psi_t\}_{0\leq t \leq 1})=-[\sigma]$$
and hence $[\sigma]\in \Gamma (M,L)$.

Conversely, let $[\sigma]\in \Gamma (M,L)$.  Choose a loop $\psi_t
\in \textmd{Symp}_0(M,L,\omega)$ such that
$\textmd{Flux}(\{\psi_t\})=-[\sigma]$. This extends to the
interval $[0,2]$ by $\mathcal{C}(\psi_t)=[(t-1)\sigma]$, for
$1\leq t \leq 2$.  Then the resulting path $\psi_t$, $t\in [0,2]$,
has zero flux.  Now by Theorem~\ref{thm-kernel} this path can be
deformed to a Hamiltonian isotopy via a homotopy with fixed end
points. Therefore $\psi=\psi_2$ is a Hamiltonian isotopy.
\end{proof}

\begin{Lemma}\label{lem-ham}
Every smooth path $\psi_t\in \textmd{\em Ham}(M,L)$ is generated
by Hamiltonian vector fields.
\end{Lemma}

The proof of the above lemma is completely analogous to that of
Proposition $10.17$ of \cite{McS} and thus will be omitted.  The
following corollary which is the relative version of Corollary
$10.18$ of \cite{McS} is the main result of this section.

\begin{Corollary}\label{cor-rel} Let $(M,\omega)$ be a closed
connected symplectic manifold and $L$ is a closed connected
Lagrangian submanifold.

{\bf i)} There is an exact sequence of simply connected Lie groups
$$1\rightarrow \widetilde{\textmd{Ham}}(M,L)\rightarrow
\widetilde{\textmd{Symp}}_0(M,L,\omega)\rightarrow
H^1(M,L,{\mathbb R})\rightarrow 0$$ where $\widetilde{\textmd{
Ham}}(M,L)$ is the universal cover of $\textmd{Ham}(M,L)$ and the
third homomorphism is the flux homomorphism.

{\bf ii)} There is an exact sequence of Lie algebras
$$0\rightarrow {\mathbb R}\rightarrow C^{\infty}(M,L)\rightarrow
\chi(M,L,\omega) \rightarrow H^1(M,L,{\mathbb R})\rightarrow 0.$$
Here the third map is $H\mapsto X_H$, the fourth map is $X\mapsto
[\imath(X)\omega]$, while $C^{\infty}(M,L)$ denotes the algebra of
smooth functions on $M$, which are constant on $L$ and
$\chi(M,L,\omega)$ is the algebra of symplectic vector fields $X$
on $M$ such that $X(q)\in T_qL$, for any $q \in L$.

{\bf iii)} The sequence of groups
$$1\rightarrow \pi_1(\textmd{Ham}(M,L))
\rightarrow \pi_1(\textmd{Symp}_0(M,L,\omega))\rightarrow \Gamma
(M,L)\rightarrow 0$$ is exact.

{\bf iv)} There is an exact sequence of groups
$$1\rightarrow \textmd{Ham}(M,L)\rightarrow
\textmd{Symp}_0(M,L,\omega)\rightarrow H^1(M,L,{\mathbb R})/
\Gamma (M,L) \rightarrow 0$$ where the third map is induced by the
flux homomorphism.
\end{Corollary}

\begin{proof}
By Lemma~\ref{lem-ham} every smooth path $\psi_t \in \textmd{
Ham}(M,L)$ which starts at the identity is a Hamiltonian isotopy
and therefore has zero flux.  This implies that
$\widetilde{\textmd{Ham}}(M,L)\subseteq \ker(\textmd{Flux})$.  On
the other hand, by Theorem~\ref{thm-kernel} if
$\textmd{Flux}(\{\psi_t\})=0$ then the path $\psi_t$ is homotopic,
with fixed end points, to a Hamiltonian isotopy and hence $\{
\psi_t\} \in \widetilde{\textmd{Ham}}(M,L)$.  Now, the first
statement follows from the surjectivity of the flux homomorphism
(see Lemma~\ref{lem-defn}).  The second statement is easy and
indeed follows from the proof of Lemma~\ref{lem-defn}.

The only nontrivial part in the third statement is the fact that
the homomorphism $\pi_1(\textmd{Ham}(M,L)) \rightarrow
\pi_1(\textmd{ Symp}_0(M,L,\omega))$ is injective.  To see this,
it suffices to show that any path $[0,1] \rightarrow
\widetilde{\textmd{Symp}}_0(M,L,\omega)$ with end points in
$\widetilde{\textmd{Ham}}(M,L)$ is isotopic with fixed end points
to a path in $\widetilde{\textmd{Ham}}(M,L)=\ker(\textmd{Flux})$.
However, this is just the parameterized version of the first
statement.

The last statement is also obvious.
\end{proof}

\section{Applications and Stability under Symplectic Reduction}

Combining Part (i) of Corollary~\ref{cor-rel} with its absolute
version
$$0\rightarrow \widetilde{\textmd{Ham}}(M)\rightarrow
\widetilde{\textmd{Symp}}_0(M,\omega)\rightarrow H^1(M,{\mathbb
R})\rightarrow 0,$$ Corollary $10.18$ of \cite{McS}, we derive the
following corollaries.

\begin{Corollary}\label{cor3}
Suppose that $(M,\omega)$ and $L\subseteq M$ are as in
Corollary~\ref{cor-rel}. If the first cohomology group,
$H^1(M,\mathbb R)$, is trivial then
$$\widetilde{\textmd{Ham}}(M)=\widetilde{\textmd{Symp}}_0(M,\omega)$$ and
$$\widetilde{\textmd{Ham}}(M,L)=\widetilde{\textmd{Symp}}_0(M,L,\omega).$$
\end{Corollary}

\begin{Corollary}\label{cor4}
Let $(M,\omega)$ be closed connected symplectic manifold and
$L\subseteq M$ a closed connected Lagrangian submanifold.  If
$ImH^1(L,{\mathbb R})$ denotes the image of the restriction map
$H^1(M,{\mathbb R})\rightarrow H^1(L,{\mathbb R})$ then the flux
homomorphism induces an exact sequence
$$1\rightarrow N(M,L)\rightarrow\widetilde{\textmd{Symp}}_0(M,\omega)
\rightarrow ImH^1(L,{\mathbb R})\rightarrow 0,$$ where  $N(M,L)$
is the normal closure of the product of the image of the canonical
map $$\widetilde{\textmd{Symp}}_0(M,L,\omega)\rightarrow
\widetilde{\textmd{Symp}}_0(M,\omega)$$  and \
$\widetilde{\textmd{Ham}}(M)$ \ in \
$\widetilde{\textmd{Symp}}_0(M,\omega)$.
\end{Corollary}

\begin{Example}\label{Exm-rel}
{\bf i)} Let $L=S^1\subseteq S^2=M$, where $S^2$ is equipped with
any symplectic form.  Then by the above corollaries
$$\widetilde{\textmd{Ham}}(M)=
\widetilde{\textmd{Symp}}_0(M,\omega),$$
$$\widetilde{\textmd{Ham}}(M,L)=
\widetilde{\textmd{Symp}}_0(M,L,\omega)$$ and
$$\widetilde{\textmd{Symp}}_0(M,\omega)=N(M,L).$$

{\bf ii)} Let $L=S^1\times \{pt\} \subseteq S^1 \times S^1=M$,
where $M$ is equipped with any symplectic form.  Then
$$\widetilde{\textmd{Symp}}_0(M,\omega)/\widetilde{\textmd{Ham}}(M)
\simeq{\mathbb R}^2,$$
$$\widetilde{\textmd{Symp}}_0(M,L,\omega)/\widetilde{\textmd{Ham}}(M,L)
\simeq{\mathbb R}$$ and
$$\widetilde{\textmd{Symp}}_0(M,\omega)/N(M,L)\simeq{\mathbb R}.$$

{\bf iii)} Again let $M=S^1\times S^1$ and $L$ any smoothly
embedded nullhomotopic circle in $M$.  Then
$$\widetilde{\textmd{Symp}}_0(M,\omega)/\widetilde{\textmd{Ham}}(M)
\simeq{\mathbb R}^2,$$
$$\widetilde{\textmd{Symp}}_0(M,L,\omega)/\widetilde{\textmd{Ham}}(M,L)
\simeq{\mathbb R}^2$$ and
$$\widetilde{\textmd{Symp}}_0(M,\omega)=N(M,L).$$
\end{Example}

\subsection{Stability under symplectic reductions}
The next result is about the stability of the quotients
$$\widetilde{\textmd{Symp}}_0(M,\omega)/\widetilde{\textmd{Ham}}(M)
\ \ \mbox{and} \ \ \widetilde{\textmd{Symp}}_0/N(M,L)$$ under
symplectic reduction.

Let $S^1$ act in a Hamiltonian fashion on the symplectic manifold
$(M,\omega)$ and $L$ be a Lagrangian submanifold contained in a
level set $\mu^{-1}(c)$ of the moment map, $\mu :M\rightarrow
{\mathbb R}$, of the $S^1$-action. Assume that $L$ is $S^1$
invariant and the $S^1$-action on the level set $\mu^{-1}(c)$ is
free.  Let $M_{red}=\mu^{-1}(c)/S^1$ be the symplectic quotient
with the symplectic structure $\omega_{red}$ and $L_{red}=L/S^1$,
which is a Lagrangian submanifold of $(M_{red},\omega_{red})$. A
theorem of Li (\cite{Li}) states that the fundamental groups and
hence the first cohomology groups of $M$, $\mu^{-1}(c)$ and
$M_{red}$ are all isomorphic under the canonical homomorphisms. On
the other hand, if $p:L\rightarrow L_{red}$ denotes the quotient
map then we have the following result (see \cite{Oz1}).

\begin{Corollary}\label{cor-final}
Let $M$, $M_{red}$, $L$ and $L_{red}$ be as above.  Then the map
$$p^*:Im(H^i(M_{red},{\mathbb Q})\rightarrow H^i(L_{red},{\mathbb Q}))
\rightarrow Im(H^i(M,{\mathbb Q})\rightarrow H^i(L,{\mathbb Q}))$$
is onto for any $i$, and is an isomorphism for $i=1$.
\end{Corollary}

The above arguments together with Corollary~\ref{cor4} imply the
following stability result.

\begin{Corollary}\label{cor-stable}
Let $M$, $M_{red}$, $L$ and $L_{red}$ be as above.  Then we have
$$\widetilde{\textmd{Symp}}_0(M,\omega)/
\widetilde{\textmd{Ham}}(M) \simeq \widetilde{\textmd{Symp}}_0
(M_{red},\omega_{red})/\widetilde{\textmd{Ham}}(M_{red})$$ and
$$\widetilde{\textmd{Symp}}_0(M,\omega)/N(M,L)\simeq
\widetilde{\textmd{Symp}}_0(M_{red},\omega_{red})/N(M_{red},L_{red}).$$
\end{Corollary}

\subsection{Real algebraic varieties and their
complexifications} In this subsection we will mention a different
kind of stability of the above groups for real algebraic
varieties. Let $X$ be a nonsingular compact real algebraic variety
with a nonsingular projective complexification $i:X \rightarrow
X_{\mathbb C}$.  Note that any two nonsingular projective
complexifications of $X$ are always birationally isomorphic and
thus have the same fundamental group and the same first cohomology
group (cf. see p.494 of \cite{GH}). Clearly $X_{\mathbb C}$
carries a K\"ahler and hence a symplectic structure such that $X$
becomes a Lagrangian submanifold. So, by Corollary $10.18$ of
\cite{McS} (or just let $L$ be the empty set in
Corollary~\ref{cor-rel}(i)) the quotient group
$$\widetilde{\textmd{Symp}}_0(X_{\mathbb C},\omega)/
\widetilde{\textmd{Ham}}(X_{\mathbb C})$$ is determined only by
$X$ and hence is independent of the projective complexification
$i:X\rightarrow X_{\mathbb C}$.

Define $KH_i(X,{\mathbb R})$ as the kernel of the induced
homomorphism
$$i_*:H_i(X,{\mathbb R}) \rightarrow H_i(X_{\mathbb C},{\mathbb
R})$$ and $ImH^i(X,{\mathbb R})$ as the image of the induced
homomorphism
$$i^*:H^i(X_{\mathbb C},{\mathbb R}) \rightarrow H^i(X,{\mathbb
R}).$$ In \cite{Oz2,Oz3} it is shown that both $KH_i(X,{\mathbb
R})$ and $ImH^i(X,{\mathbb R})$ are independent of the projective
complexification $i:X \rightarrow X_{\mathbb C}$ and thus (entire
rational) isomorphism invariants of $X$.  Moreover, by
Corollary~\ref{cor4}, for any topological component $X_0$ of $X$
the quotient group
$$\widetilde{\textmd{Symp}}_0(X_{\mathbb C},\omega)/N(X_{\mathbb
C},X_0)$$ is independent of the smooth projective complexification
$i:X\rightarrow X_{\mathbb C}$ and hence is determined only by $X$
or equally by $X_0$.

\bibliographystyle{amsplain}

\end{document}